\documentclass[11pt]{article}
\usepackage[english]{babel}
\usepackage{amssymb, amsmath, amsthm} 
\usepackage[a4paper, centering]{geometry}
\usepackage{graphicx,color}
\usepackage{subfigure}

\usepackage{amsthm}
\makeatletter
\def\th@plain{%
  \thm@notefont{}
  \itshape 
}
\def\th@definition{%
  \thm@notefont{}
  \normalfont 
}
\makeatother

\newtheorem{thm}{Theorem}[section]
\newtheorem{prop}[thm]{Proposition}

\theoremstyle{definition}
\newtheorem{defn}[thm]{Definition}
\newtheorem{example}[thm]{Example}
\newtheorem{oss}[thm]{Remark}

\newcommand{\Z}{\mathbb{Z}}
\newcommand{\N}{\mathbb{N}}
\newcommand{\R}{\mathbb{R}}

\renewcommand{\H}{\mathcal{H}}
\newcommand{\I}{\mathcal{I}}
\newcommand{\D}{\mathcal{D}}
\newcommand{\Na}{N_\alpha}
\newcommand{\Nb}{N_\beta}
\renewcommand{\P}{\text{P}}
\renewcommand{\epsilon}{\varepsilon}

\def\e{\epsilon}

\begin{document}

\title{Motion of discrete interfaces in periodic media}

\author{{\scshape Andrea Braides }\\
Dipartimento di Matematica, Universit\`a di Roma `Tor Vergata'\\
via della ricerca scientifica 1, 00133 Roma (Italy)\\ \\
{\scshape Giovanni Scilla}\\
Dipartimento di Matematica `G.~Castelnuovo'\\
`Sapienza' Universit\`a di Roma\\
piazzale Aldo Moro 5, 00185 Roma (Italy)}
\date{}

\maketitle

\begin{abstract}\noindent
We study the motion of discrete interfaces driven by ferromagnetic interactions in a two-dimensional periodic environment by coupling the minimizing movements approach by Almgren, Taylor and Wang and a discrete-to-continuous analysis. 
The case of a homogeneous environment has been recently treated by Braides, Gelli and Novaga, showing that the effective continuous motion is a flat motion related to the crystalline perimeter obtained by $\Gamma$-convergence
from the ferromagnetic energies, with an additional discontinuous dependence on the curvature, giving in particular a pinning threshold. In this paper we give an example showing that in general the motion does not depend only on the $\Gamma$-limit, but also on geometrical features that are not detected in the static description. In particular we show 
how the pinning threshold is influenced by the microstructure and that the effective motion is described by a new homogenized velocity.
\end{abstract}

\section{Introduction}
In this paper we study a model problem of homogenization for a discrete crystalline flow. The analysis will be
carried over by using the minimizing-movement scheme of Almgren, Taylor and Wang \cite{ATW83} 
(later thus renamed by De Giorgi, see e.g.~\cite{AGS}). This consists in introducing a time scale $\tau$, iteratively defining
a sequence of sets $E^\tau_k$ as minimizers of  
\begin{equation}
\min \Bigl\{ P(E)+{1\over \tau} D(E, E^\tau_{k-1})\Bigr\} ,
\end{equation}
where $P$ is a perimeter energy and $D$ is a (suitably defined) distance-type energy between sets, and $E^\tau_0$ 
is a given initial datum, and subsequently computing a time-continuous limit $E(t)$ of $\{E^\tau_k\}$ as $\tau\to0$, which defines the desired geometric motion related to the energy~$P$. 

The study of geometric motions in inhomogeneous 
environments has a very large literature (see e.g.~\cite{BC, DY, Gl, LS}). The ones in a discrete setting can be 
considered a somewhat extreme case, in that the corresponding energies possess a large number of local minimizers 
(actually, by the discrete nature of the problem all states are local minimizers), while on the contrary their continuum limits
(see e.g.~\cite{ABC,BraPia2} for a rigorous definition) possess no local (non global) minimizer. As a consequence, gradient flows tend to be ``pinned'' 
(i.e., the resulting limit $E(t)$ is constant), in contrast with the formal evolution of their limit continuous energies, 
to which the Almgren-Taylor-Wang approach can be used to obtain a non trivial evolution (for the case of $P$ the crystalline perimeter in two dimensions see \cite{AT95}).

In a recent paper by Braides, Gelli and Novaga \cite{BGN} the Almgren-Taylor-Wang approach has been used coupled to a homogenization procedure. In this case the perimeters (and the distances) depend on a small parameter $\e$, and consequently, after introducing a time scale $\tau$, the time-discrete motions are the $E^{\tau,\e}_k$ defined iteratively by 
\begin{equation}
E^{\tau,\e}_k \hbox{ is a minimizer of } \min \Bigl\{ P_\e(E)+{1\over \tau} D(E, E^{\tau,\e}_{k-1})\Bigr\}.
\end{equation}
The time-continuous limit $E(t)$ of $\{E^{\tau,\e}_k\}$ then may depend how mutually $\e$ and $\tau$ tend to $0$.
This type of problems can be cast in the general framework of {\em minimizing movements along a $\Gamma$-converging sequence}  (see \cite{Bra13}). In particular, if we have a large number of local minimizers then the limit motion will be pinned if $\tau<\!\!<\e$ suitably fast (in a sense, we can pass to the limit in $\tau$ first, and then apply the Almgren-Taylor-Wang approach, which clearly gives pinning when the initial data are local minimizers). On the contrary, if $\e<\!\!<\tau$ fast enough 
and $P_\e$ $\Gamma$-converge to a limit perimeter $P$ (which is always the case by compactness), then the limit $E$ will be 
the evolution related to the limit $P$ (again, in a sense, in this case we can pass to the limit in $\e$ first). 

In \cite{BGN} the energies $P_\e$ are {\em ferromagnetic energies} defined on subsets $E\subset \e\Z^2$, of the form
$$
P_\e(E)= {\alpha}\,\e\, \#\{(i,j)\in\e\Z^2\times\e\Z^2: i\in E, j\not\in E, \ |i-j|=\e\}
$$
($\alpha>0$ a positive parameter). The continuum limit of these energies can be proved to be the crystalline perimeter
$$
P(E)=\alpha\int_{\partial E}\|\nu\|_1d\H^1,
$$
where $\nu$ is the normal to $\partial E$ and $\|(\nu_1,\nu_2)\|_1=|\nu_1|+|\nu_2|$ (see \cite{ABC}). The {\em flat flow} of this perimeter is the motion by crystalline curvature described by Taylor \cite{Ta}. In the case of initial datum a coordinate rectangle, the evolution is a rectangle with the same centre and sides of lengths $L_1, L_2$ governed by the system of ordinary differential equations
$$
\begin{cases}\displaystyle \dot L_1= -{4\alpha\over L_2}\cr\cr
\displaystyle \dot L_2= -{4\alpha \over L_1}.\end{cases}
$$

In \cite{BGN} all possible evolutions have been characterized as $\e, \tau\to 0$, showing that the relevant mutual scale is when $\tau/\e\to\gamma\in(0,+\infty)$. In the case of initial datum a coordinate rectangle the resulting evolution is still a rectangle. In the case of a unique evolution, the side-lengths $L_1(t), L_2(t)$ of this rectangle are governed by a system of `degenerate' ordinary differential equations
$$
\begin{cases}\displaystyle \dot L_1= -{2\over \gamma}\Bigl\lfloor{2\gamma\alpha\over L_2}\Bigr\rfloor\\ \cr
\displaystyle \dot L_2= -{2\over \gamma}\Bigl\lfloor{2\gamma\alpha\over L_1}\Bigr\rfloor.\end{cases}
$$
Note that the right-hand sides are discontinuous; however existence (and uniqueness, except for some special cases) of solution can be easily proved by a direct argument. This characterization highlights the effect of the periodicity through the scaling $\gamma$ and that the motion is obtained by overcoming some energy barriers in a `quantized' manner by the presence of a discontinuous right-hand side. In particular, we have {\em pinning of large rectangles}: if both initial side-lengths are above the {\em pinning threshold} $\widetilde L=2\gamma\alpha$ then the right-hand sides are zero and the motion is pinned. The limit cases (total pinning and continuous crystalline flow) correspond to the limit values $\gamma=0$ and $\gamma=+\infty$. This analysis shows that the ``correct scaling'' for this problem is $\e\sim\tau$, which gives the most information about all the limit evolutions.

The analysis described above exhibits a limit evolution in which we may read the effect of the $\Gamma$-limit energy (through the crystalline form of the evolution and the coefficient $\alpha$) and of the interplay between the time and space scales through the scaling $\gamma$. Scope of this work is to show that in general the situation can be more complex, and the periodic microstructure can affect the limit evolution without changing the $\Gamma$-limit. To this end we will introduce a further inhomogeneity in the perimeters $P_\e$ by considering
$$
P_\e(E)= {1\over 2}\e\,\sum \{c_{ij}:i,j\in \Z^2,\e i\in E, \e j\not\in E, \ |i-j|=1\},
$$
(we use the notation $\sum \{ x_a: a\in A\}=\sum_{a\in A} x_a$)
where the coefficients $c_{ij}$ equal $\alpha$ except for some well-separated periodic square 
inclusions where $c_{ij}=\beta>\alpha$.
These inclusions are not energetically favorable and they can be neglected in the computation of the $\Gamma$-limit, which is still the perimeter $P$ above, with the same coefficient $\alpha$. They can be considered as ``obstacles'' that can be bypassed when computing minimizers of $P_\e$; however their presence is felt in the minimizing-movement procedure since they may influence the choice of $E^{\tau,\e}_k$ through the interplay between the distance and perimeter terms. As a result, the motion can be either decelerated or accelerated with respect to the homogeneous case.

As already remarked in \cite{BGN} the relevant case for the description of the motion is that of initial data coordinate 
rectangles, since all other cases can be reduced to the study of this one. We will then restrict our analysis to that case. 
This (apparently) simple situation already contains all the relevant features of the evolution and highlights the differences with respect to \cite{BGN}. We will show that the limit motion can still be described through a system of degenerate ordinary differential equations of the form
$$
\begin{cases}\displaystyle \dot L_1= -{2\over \gamma}\,f\Bigl({\gamma\over L_2}\Bigr)\cr \cr
\displaystyle \dot L_2= -{2\over \gamma}\,f\Bigl({\gamma\over L_1}\Bigr)\end{cases}
$$
with $f$ a locally constant function on compact subsets of $(0,+\infty)$ which depends on $\alpha$, the period and size of the inclusions but not on $\gamma$ (neither on the value $\beta$). The {\em effective velocity} $f$ is obtained by a homogenization formula which optimizes the motion of the sides of the rectangle, resulting in an oscillation around a linear motion with velocity ${1\over\gamma} f(\gamma/L)$ (which is locally constant as noted above). Note that, in the case of no inclusion, the system is of the same form with $f(Y)=\lfloor 2\alpha Y\rfloor$.
The dependence on the inclusions gives a new pinning threshold
$$
\overline L={4\gamma\alpha\over  2+\Nb}
$$
depending on the size of the inclusion $\Nb$.
The reason for this new pinning threshold is that, in order that a side may move, it needs to be able to overcome a barrier of $\Nb$ inclusions. Note that, if the initial data have side-lengths $\overline L< L<\widetilde L$, then we may have a microscopic motion which stops after a finite number of time steps, and is not eventually detected in the limit. It should be remarked that the presence of the inclusions may indeed accelerate the motion, so that $f(Y)>\lfloor 2\alpha Y\rfloor$ for some $Y$.

\bigskip

The paper is organized as follows. In Section~\ref{setting} we define all the energies that we will consider. We then formulate the discrete-in-time scheme analogous to the Almgren, Taylor and Wang approach. Section~\ref{rectangular} contains the proof of the convergence of the discrete scheme in the case of a rectangular initial set. Contrary to the case in \cite{BGN} it is not trivial to show that the minimizers of this scheme are actually rectangles. This is a technical result contained in Proposition~\ref{rectangleprop}. Subsection~\ref{newpinning} contains the computation of the new pinning threshold, showing that it depends on the percentage $\Nb$ of defects in the lattice. Subsection~\ref{newvelocity} deals with the new definition of the effective velocity of a side by means of a homogenization formula resulting from a one-dimensional `oscillation-optimization' problem. This velocity can be expressed uniquely (up possibly to a discrete set of values), as a function the ratio of $\gamma$ and the side-length (Definition~\ref{effvel}). The description of the homogenized limit motion is contained in Subsection~\ref{limitmotion}. In the last Section 4 we explicitly compute the velocity function by means of algebraic formulas in some simple cases, showing a nontrivial comparison with the case with no inclusions.

\section{Setting of the problem}\label{setting}

If $x=(x_1,x_2)\in \R^2$ we set $\|x\|_1=|x_1|+|x_2|$ and $\|x\|_\infty=\max\{|x_1|,|x_2|\}$. If $A$ is a Le\-bes\-gue\hbox{-}measurable set we denote by $|A|$ its two\hbox{-}dimensional Lebesgue measure. The symmetric difference of $A$ and $B$ is denoted by $A\triangle B$, their Hausdorff distance by $\text{d}_\mathcal{H}(A,B)$.
If $E$ is a set of finite perimeter then $\partial^*E$ is its reduced boundary (see, for example \cite{Bra98}). The 
measure-theoretical inner normal to $E$ at a point $x$ in $\partial^*E$ is denoted by $\nu=\nu_E(x)$ .

\subsection{Inhomogeneous ferromagnetic  energies}
\begin{figure}[htbp]
\centering
\def\svgwidth{100pt}
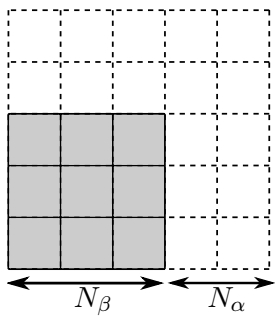
\caption{Periodicity cell. Continuous lines represent $\beta$\hbox{-}bonds, dashed lines $\alpha$\hbox{-}bonds.}\label{fig:0}
\end{figure}

\def\Nab{N_{\alpha\beta}}
The energies we consider are interfacial energies defined in an inhomogeneous environment as follows: let $0<\alpha<\beta<+\infty$, $\Na,\Nb\geq1$ and set $N_{\alpha\beta}=\Na+\Nb$. We consider the $\Nab$-periodic coefficients $c_{ij}$ indexed on {\em nearest-neighbours} of $\Z^2$ (i.e., $i,j\in \Z^2$  with $|i-j|=1$) defined for $i,j$ such that 
$$
0\le  {i_1+j_1\over 2}, {i_2+j_2\over 2}<\Nab
$$
by 
\begin{equation}
c_{ij}= \begin{cases} \beta &\hbox{if $\displaystyle0\le {i_1+j_1\over 2}, {i_2+j_2\over 2}\le \Nb$}\\ 
\alpha &\hbox{otherwise.}\end{cases}
\end{equation}
These coefficients label the bonds between points in $\Z^2$, so that they describe a matrix of $\alpha$-bonds with $\Nab$-periodic inclusions of $\beta$-bonds grouped in squares of side-length $\Nb$. The periodicity cell is pictured in Fig.~\ref{fig:0}.
Correspondingly, to these coefficients we associate the energy defined on subsets $\I$ of $\Z^2$ by
\begin{equation}
\P^{\alpha,\beta}(\I)=\sum\Bigl\{c_{ij}: |i-j|=1, i\in \I, j\in \Z^2\setminus \I\Bigr\}.
\end{equation}
As recalled in the Introduction we use the notation $\sum \{ x_a: a\in A\}=\sum_{a\in A} x_a$.

In order to examine the overall properties of $\P^{\alpha,\beta}$ we introduce the family of scaled energies defined on subsets $\I$ of $\e\Z^2$ by
\begin{equation}\label{pabe}
\P^{\alpha,\beta}_\epsilon(\I)=\sum\Bigl\{\e\ c_{{i/\e}\,{ j/\e}}: |i-j|=\e, i\in \I, j\in \e\Z^2\setminus \I\Bigr\};
\end{equation}
i.e., $\P^{\alpha,\beta}_\e(\I)=\e\, \P^{\alpha,\beta}({1\over\e}\I)$.
To study the continuous limit as $\e\to0$ of these energies it will be convenient to identify each subset of $\e\Z^2$ with a measurable subset of $\R^2$, in such a way that equi-boundedness of the energies implies pre-compactness of such sets in the sense of sets of finite perimeter. This identification is as follows: we denote by $Q$ the closed coordinate unit square of center $0$, $Q=[-1/2,1/2]^2$; if $\epsilon>0$ and $ {i}\in\epsilon\Z^2$, we denote by $Q_\epsilon( {i})= {i}+\epsilon Q$ the closed coordinate square with side-length $\epsilon$ and centered in $ {i}$. To a set of indices $\I\subset\epsilon\Z^2$ we associate the set
\begin{equation*}
E_{\I}=\bigcup_{ {i}\in \I}Q_\epsilon( {i}).
\end{equation*} 
\\
The space of \emph{admissible sets} related to indices in the two\hbox{-}dimensional square lattice is then defined by
\begin{equation*}
\D_\epsilon:=\left\{E\subseteq\R^2:\quad E=E_\I\text{ for some $\I\subseteq\epsilon\Z^2$}\right\}.
\end{equation*}
For each $E=E_\I\in \D_\epsilon$ we denote
\begin{equation}\label{pabei}
\P^{\alpha,\beta}_\epsilon(E)=\P^{\alpha,\beta}_\epsilon(\I).
\end{equation}
As an easy remark, we note that 
\begin{equation}\label{pabes}
\P^{\alpha,\beta}_\epsilon(E)\ge\e\alpha\#\Bigl\{(i,j): |i-j|=\e, i\in \I, j\in \e\Z^2\setminus \I\Bigr\}= \alpha \H^1(\partial E),
\end{equation}
which shows that sequences of sets $E_\e$ with $\sup_\e\P^{\alpha,\beta}_\epsilon(E_\e)<+\infty$ are pre-compact with respect 
to the local $L^1$-convergence in $\R^2$ of their characteristic function and their limits are sets of finite perimeter in $\R^2$.
Hence, this defines a meaningful convergence with respect to which compute the $\Gamma$-limit of $\P^{\alpha,\beta}_\epsilon$ as $\e\to 0$.

A general theory for the homogenization of energies (\ref{pabe}), in a more general context, has been developed in \cite{BraPia2} (see also \cite{BCS,BraPia3,BraSo1}), where it is shown that the $\Gamma$-limit's domain is precisely the family of sets of finite perimeter and its general form is
$$
F(E)=\int_{\partial^*E}\varphi(\nu)d\H^1,
$$
with $\varphi$ a convex function positively homogeneous of degree one.
The computation in the case $\alpha=\beta$ (homogeneous spin systems) can be found in \cite{ABC} and gives
$\varphi(\nu)=\alpha\|\nu\|_1$. In our case the presence of the $\beta$-inclusions does not influence the form of the $\Gamma$-limit, as in the following remark.

\begin{oss} [$\Gamma$\hbox{-}convergence of inhomogeneous perimeter energies]
The energies $\text{P}_\epsilon^{\alpha,\beta}$ defined by (\ref{pabe}) $\Gamma$\hbox{-}converge, as $\epsilon\to0$, to the anisotropic crystalline perimeter functional
\begin{equation*}
\P^\alpha(E)=\alpha\int_{\partial^* E}\|\nu\|_1\,d\H^1.
\end{equation*}
This limit is independent of $\Na,\Nb$, and equals the one obtained when $\beta=\alpha$. 

The lower bound for the $\Gamma$-limit is immediately obtained from the case $\alpha=\beta$ in \cite{ABC} after remarking that $\P^{\alpha,\beta}_\epsilon\geq \P_\epsilon^{\alpha,\alpha}$. In order to verify the upper bound, it suffices to note that recovery sequences for the $\Gamma$-limit of $\P_\epsilon^{\alpha,\alpha}$ can be constructed at a scale $\Nab\e$, thus `avoiding' the $\beta$-connections. To this end, define 
$$
{Q}^{\Nab}_\epsilon=\displaystyle\bigcup\Bigl\{Q_\epsilon( {i}):\ 
{ {i}\in\epsilon\Z^2},\ 0\leq\| {i}\|_{\infty}<\e N_{\alpha\beta}\Bigr\}.
$$
This is a square of side-length $\Nab\,\e$ whose boundary intersects only $\alpha$-bonds.
We consider $\P^{\Nab}_\epsilon$the restriction of $\P_\epsilon^{\alpha,\beta}$ to the class
\begin{equation*}
{\D}^{\Nab}_\epsilon=\left\{E\subseteq\R^2:\quad E\text{ is a finite union of  $\e\Z^2$-translations of }{Q}^{\Nab}_\epsilon\right\}.
\end{equation*}
Note that we have $\P_\epsilon^{\alpha,\beta}(E)=\P_\epsilon^{\alpha,\alpha}(E)$ for $E\in {\D}^{\Nab}_\epsilon$, and that
sets in ${\D}^{\Nab}_\epsilon$ differ from sets in ${\D}_{\e\Nab}$ by a fixed translation of order  $\e$. Hence, we have
(see \cite{GCB} for details on the properties of $\Gamma$-upper limits)
$$
\Gamma\hbox{-}\limsup_{\e\to 0}\P_\epsilon^{\alpha,\beta}(E)\le 
\Gamma\hbox{-}\limsup_{\e\to 0}\P^{\Nab}_\epsilon(E)=\Gamma\hbox{-}\lim_{\e\to 0}\P_{\Nab\epsilon}^{\alpha,\beta}(E),
$$
and the latter is again equal to $\P^\alpha(E)$. This inequality just states that we can take sets in ${\D}^{\Nab}_\epsilon$ which are (small translations of) a recovery sequence for $\P_{\Nab\epsilon}^{\alpha,\beta}(E)$ as  a recovery sequence for $\P_\epsilon^{\alpha,\beta}(E)$.
\end{oss}

\subsection{A discrete-in-time minimization scheme}\label{timemin} 

For $\I\subset\epsilon\Z^2$ we define the \emph{discrete $\ell^{\infty}$\hbox{-}distance} from $\partial\I$ as

\begin{equation*}
d_\infty^\epsilon( {i},\partial\I)=
\begin{cases}
\inf\{\| {i}- {j}\|_\infty: {j}\in\I\}&\text{if $ {i}\not\in\I$}\\
\inf\{\| {i}- {j}\|_\infty: {j}\in\e\Z^2\backslash\I\}&\text{if $ {i}\in\I$}.
\end{cases}
\end{equation*}
Note that we have
$\displaystyle 
d_\infty^\epsilon( {i},\partial\I)=d_\infty( {i},\partial E_\I)+\frac{\epsilon}{2}$,
where $d_\infty$ denotes the usual $\ell^\infty$-distance. The distance can be extended to all $\R^2\backslash\partial E_\I$ by setting

\begin{equation*}
d_\infty^\epsilon(x,\partial\I)=d_\infty^\epsilon( {i},\partial\I)\quad \text{if }x\in Q_\epsilon( {i}).
\end{equation*}
\\
In the following we will directly work with $E\in\D_\epsilon$, so that the distance can be equivalently defined by

\begin{equation*}
d_\infty^\epsilon(x,\partial E)=d_\infty( {i},\partial E)+\frac{\epsilon}{2},\quad \text{if }x\in Q_\epsilon( {i}).
\end{equation*}
Note that this is well defined as a measurable function, since its definition is unique outside the union of the boundaries of the squares $Q_\epsilon$ (that are a negligible set).\\

We now fix a time step $\tau>0$ and introduce a discrete motion with underlying time step $\tau$ obtained by successive minimization. At each time step we will minimize an energy $\mathcal{F}_{\epsilon,\tau}^{\alpha,\beta}:\D_\epsilon\times\D_\epsilon\to\R$ defined as
\begin{equation}
\mathcal{F}_{\epsilon,\tau}^{\alpha,\beta}(E,F)= \P_\epsilon^{\alpha,\beta}(E)+\frac{1}{\tau}\int_{E\triangle F}d_\infty^\epsilon(x,\partial F)\,dx.
\label{newenergy}
\end{equation}
Note that the integral can be indeed rewritten as a sum on the set of indices $\e\Z^2\cap(E\triangle F)$ (see \cite{BGN}).

Given an initial set $E^0_{\epsilon}$, we define recursively a sequence $E_{\epsilon,\tau}^k$ in $\D_\epsilon$ by requiring the following:
\begin{description}
\item[(i)] $E^0_{\epsilon,\tau}=E^0_{\epsilon}$;
\item[(ii)] $E_{\epsilon,\tau}^{k+1}$ is a minimizer of the functional $\mathcal{F}_{\epsilon,\tau}^{\alpha,\beta}(\cdot,E_{\epsilon,\tau}^k)$.
\end{description}
The \emph{discrete flat flow} associated to functionals $\mathcal{F}_{\epsilon,\tau}^{\alpha,\beta}$ is thus defined by

\begin{equation}\label{disefo}
E_{\epsilon,\tau}(t)=E_{\epsilon,\tau}^{\lfloor t\slash\tau\rfloor}.
\end{equation}
Assuming that the initial data $E^0_\epsilon$ tend, for instance in the Hausdorff sense, to a sufficiently regular set $E_0$, we are interested in identifying the motion described by any converging subsequence of $E_{\epsilon,\tau}(t)$ as $\epsilon,\tau\to0$. 

As remarked in the Introduction, the interaction between the two discretization parameters, in time and space, plays a relevant role in such a limiting process. More precisely, the limit motion depends strongly on their relative decrease rate to 0. If $\epsilon\!<\!<\tau$, then we may first let $\epsilon\to0$, so that $\P_\epsilon^{\alpha,\beta}(E)$ can be directly replaced by the limit anisotropic perimeter $\P^\alpha(E)$ and $\frac{1}{\tau}\int_{E\bigtriangleup F}d_\infty^\epsilon(x,\partial F)\,dx$ by $\frac{1}{\tau}\int_{E\bigtriangleup F}d_\infty(x,\partial F)\,dx$. As a consequence the approximated flat motions tend to the solution of the continuous ones studied by {Almgren} and {Taylor} \cite{AT95}. On the other hand, if $\epsilon\!>\!>\tau$ then there is no motion and $E_{\epsilon,\tau}^k\equiv E^0_\e$. Indeed, for any $F\neq E^0_\e$ and for $\tau$ small enough we have

\begin{equation*}
\frac{1}{\tau}\int_{E^0_\e\bigtriangleup F}d_\infty^\epsilon(x,\partial F)\,dx\geq c\frac{\epsilon}{\tau}>\P_\epsilon^{\alpha,\beta}(E^0_\e).
\end{equation*}
In this case the limit motion is the constant state $E_0$. The meaningful regime is the intermediate case $\tau\sim\epsilon$.

\bigskip

\section{Motion of a rectangle}\label{rectangular}

As shown in \cite{BGN} the relevant case is when $\e$ and $\tau$ are of the same order and the initial data are coordinate rectangles $E^0_\epsilon$, which will be the content of this section.

We assume that
\begin{equation*}
\tau=\gamma\epsilon\quad\text{for some }\gamma\in(0,+\infty),
\end{equation*}
and, correspondingly, we omit the dependence on $\tau$ in the notation of 
\begin{equation*}
E^k_\epsilon=E^k_{\epsilon,\tau}(=E^k_{\epsilon,\gamma\epsilon}).
\end{equation*}

\smallskip
Due to the lack of uniqueness of minimizers in the discrete minimization scheme, a standard comparison principle cannot hold. We recall a \emph{weak comparison principle} for our motion in the discrete case (see \cite{BGN} for the proof).

\begin{prop}[Discrete weak comparison principle]\label{wcp}
Let $\epsilon>0$ and let $R_\epsilon,K_\epsilon\in\D_\epsilon$ be such that $R_\epsilon\subseteq K_\epsilon$ and $R_\epsilon$ is a coordinate rectangle. Let $K_\epsilon^k$ be a motion from $K_\epsilon$ constructed by successive minimizations. Then $R_\epsilon^k\subseteq K_\epsilon^k$ for all $k\geq1$, where $R_\epsilon^k$ is a motion from $R_\epsilon$ constructed by successively choosing a minimizer of $\mathcal{F}_{\epsilon,\tau}^{\alpha,\alpha}(\cdot,R_\epsilon^{k-1})$ having smallest measure.
\end{prop}

\begin{oss}\label{wpcoss}
The set $\R^2\backslash K_\epsilon^k$ is the $k$-step evolution of the complementary $\R^2\backslash K_\epsilon^k$ of $K_\epsilon$. As a consequence, if we have $R_\epsilon\subseteq\R^2\backslash K_\epsilon^k$, from Proposition \ref{wcp} it follows that

\begin{equation*}
R_\epsilon^k\subseteq\R^2\backslash K_\epsilon^k,\quad \text{for all $k\geq1$.}
\end{equation*}
\end{oss}

\begin{defn}[$\alpha$-type rectangle] 
A coordinate rectangle whose sides intersect only $\alpha$\hbox{-}bonds will be called an \emph{$\alpha$-type rectangle}.
\end{defn}

The first result is that coordinate rectangles evolve into $\alpha$-type 
rectangles.

\begin{prop}\label{rectangleprop}
If $E^0_{\epsilon}\in\D_\epsilon$ is a coordinate rectangle and $F$ is a minimizer for the minimum problem for $\mathcal{F}^{\alpha,\beta}_{\epsilon,\tau}(\cdot,E^k_\epsilon), k\geq0$, then for all $\delta>0$ $F$ is a coordinate $\alpha$-type rectangle as long as the sides of $E^k_\epsilon$ are larger than $\delta$ and $\e$ is small enough. 
\end{prop}

\proof
{\bf Step 1: connectedness of $F$.} We want to prove that each $E_\epsilon^k$ is connected. It will suffice to show this for $F=E_\epsilon^1$. We first need an estimate on the area of the ``small components'' of $E_\epsilon^1$; this estimate will be obtained by using the comparison principle in Proposition \ref{wcp}.

\begin{figure}[ht]
\centering
\def\svgwidth{190pt}
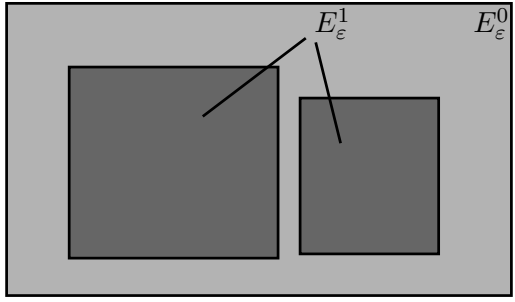
\caption{Test set with  $E^1_\epsilon$ with two large components.}\label{fig:1}
\end{figure}

Let $\ell>0$ be the maximum number such that for each point $x\in E_\epsilon^0$ there exists $y\in \R^2$ such that $x\in (y+ Q_\ell)\subseteq E_\epsilon^0$, where $Q_\ell=[-\ell/2,\ell/2]\times[-\ell/2,\ell/2]$, and the same property holds for $x\not\in E_\epsilon^0$. If $E_\epsilon^0=[-L_1/2,L_1/2]\times[-L_2/2,L_2/2]$, we can choose $\ell=\min\{L_1,L_2\}$. By applying Proposition \ref{wcp} and Remark \ref{wpcoss} to the union of squares contained in $E_\epsilon^0$, and to those outside $E_\epsilon^0$, respectively, and taking into account that a side of length $\ell$ shrinks by $\left\lfloor\frac{2\alpha\gamma}{\ell}\right\rfloor\epsilon$ in absence of defects (see \cite{BGN}), it follows that

\begin{equation*}
\text{d}_\H(\partial E_\epsilon^1,\partial E_\epsilon^0)\leq\left(\frac{2\alpha\gamma}{\ell}+1\right)\epsilon.
\end{equation*}
\\
In this way, it is not possible to have a configuration as in Fig.~\ref{fig:1}, with two large components for $E^1_\epsilon$.

\begin{figure}[ht]
\centering
\def\svgwidth{190pt}
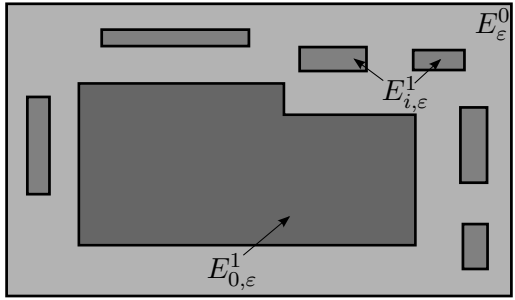
\caption{Small components of $E_\epsilon^1$.}\label{fig:2}
\end{figure}

Assume by contradiction that $E_\epsilon^1$ is not connected. In this case we should have
only one large component as in Fig.~\ref{fig:2}. 
We consider the decomposition 
\begin{equation*}
E_\epsilon^1=E_{0,\epsilon}^1\cup\bigcup_{i=1}^NE_{i,\epsilon}^1,
\end{equation*}
with $E_{0,\epsilon}^1$ the component containing all the points of $E_\epsilon^0$ having distance more than $C'\epsilon$ from $\partial E_\epsilon^0$ for a suitable constant $C'<2\alpha\gamma/\ell+1$.

Therefore for a suitable constant $C''$ we have
\begin{equation*}
d_\infty^\epsilon(x,\partial E_\epsilon^0)\leq C''\epsilon \quad\text{ for all $x\in E_{i,\epsilon}^1$ and $i\geq1$}.
\end{equation*}
By using the isoperimetric inequality, for $\epsilon$ small enough we infer
\begin{equation*}
\frac{1}{\tau}\int_{E_{i,\epsilon}^1}d_\infty^\epsilon(x,\partial E_\epsilon^0)\,dx\leq(C''/\gamma)|E_{i,\epsilon}^1|<C_{\text{iso}}\sqrt{|E_{i,\epsilon}^1|}\leq \P^{\alpha,\alpha}_\epsilon(E_{i,\epsilon}^1)\leq\P_\epsilon^{\alpha,\beta}(E^1_{i,\epsilon}),
\end{equation*}
with $C_{\text{iso}}$ being the constant of the isoperimetric inequality. Thus, we get a contradiction since we can decrease strictly the energy by eliminating the small components of $E_\epsilon^1$ and considering the set $E'=E_{0,\epsilon}^1$ as a competitor.

\begin{figure}[htbp]
\centering
\def\svgwidth{190pt}
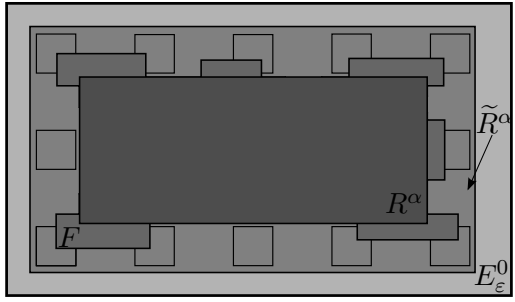
\caption{$\alpha$-rectangularization.}\label{fig:4}
\end{figure}

\smallskip
{\bf Step 2: $\alpha$\hbox{-}rectangularization.} Consider the maximal $\alpha$-type rectangle $R^\alpha$ with each side intersecting $F$. We call the set $F'=F\cup R^\alpha$ the \emph{$\alpha$\hbox{-}rectangularization} of $F$. This set is either an $\alpha$-type rectangle (and in this case we conclude) or it has some protrusions intersecting $\beta$-bonds (Fig.~\ref{fig:4}). 
In both cases $\P^{\alpha,\beta}_\epsilon(F')\leq\P^{\alpha,\beta}_\epsilon(F)$, and the symmetric difference with $E_\epsilon^0$ decreases. To justify this, note that the $\alpha$\hbox{-}rectangularization reduces (or leaves unchanged) $\P^{\alpha,\alpha}_\epsilon$ and it reduces the symmetric difference.

As a consequence of this observation, we also deduce an \emph{a priori} estimate on the maximal distance between $\partial E^0_\e$ and $\partial E^1_\e$. By the argument above, $F$ contains an $\alpha$-type rectangle $R^\alpha$
and is strictly contained in an $\alpha$-type rectangle $\widetilde R^\alpha$ whose sides have a distance from the corresponding sides of $R^\alpha$ of not more than $(\Nb+1)\e$.
We only check the \emph{a priori} estimate in the simplifying hypothesis that  $E^0_\e$ is of $\alpha$-type and that $E^0_\e$ and $R^\alpha$ are both concentric squares, so that we can express this estimate in terms of the length $L$ of the sides of $E^0_\e$ and the distance between $\partial E^0_\e$ and $\partial  R^\alpha$, which can be expressed as $\e N$. Note that we have 
$$
\alpha \H^1(\partial E^0_\e)\ge \mathcal{F}_{\epsilon,\tau}^{\alpha,\beta}(E^1_\e,E^0_\e)\ge 
\alpha \H^1(\partial R^\alpha)+\frac{1}{\tau}\int_{E^0_\e\setminus \widetilde{R}^\alpha}d_\infty^\epsilon(x,\partial E^0_\e)\,dx,
$$ 
which translates into
$$
4\alpha L\ge 4\alpha(L-2\e N)+{2L\over\gamma}\e(N-N_\beta)^2+ O(\e^2),
$$
and gives (for $\e$ sufficiently small)
\begin{equation}\label{cielle}
N\le {c_1\over L}+c_2N_\beta=:c(L).
\end{equation}
The same type of estimate holds in the general case taking $L$ the minimal length of sides of $E^0_\e$.

\smallskip
{\bf Step 3:  profile of protrusions on $\beta$\hbox{-}squares.}
Now we want to describe the form of the optimal profiles of the boundary of $F$ intersecting $\beta$-squares.

As noted above, $F$ contains an $\alpha$-type rectangle $R^\alpha=[\e m_1,\e M_1]\times [\e m_2,\e M_2]$ and is contained in the $\alpha$-type rectangle $$\widetilde R^\alpha=[\e (m_1-N_\beta),\e (M_1+N_\beta)]\times[\e (m_2-N_\beta),\e (M_2+N_\beta)]$$ whose side exceed the ones of $R^\alpha$ by at most $2\e N_\beta$. We will describe separately the possible profile of $F$ close to each side of $R^\alpha$; e.g., in the rectangle $[\e (m_1-N_\beta),\e (M_1+N_\beta)]\times[\e M_2,\e (M_2+N_\beta)]$ (i.e., close to the upper horizontal side of $R^\alpha$).

\begin{figure}[ht]
\centering
\def\svgwidth{280pt}
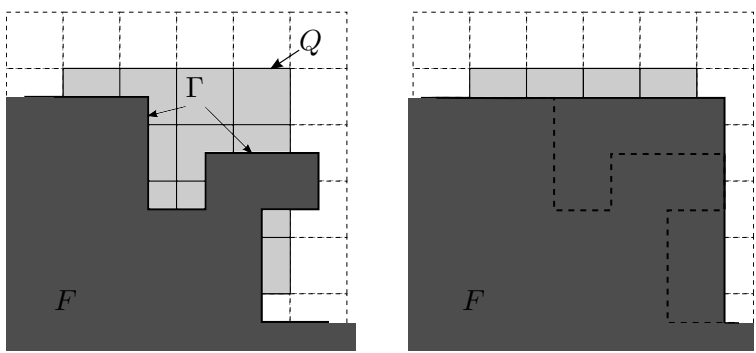
\caption{Envelope of $\partial F$ when intersecting opposite sides.}\label{upper}
\end{figure}

\bigskip
We first consider the possible behavior of the boundary of $F$ at a single $\beta$-square~$Q$. 
We suppose that such $Q$ is not one of the two extremal squares, for which a slightly different analysis holds.
First, if a portion $\Gamma$ of $\partial F$ intersects $Q$ in exactly two points on opposite vertical sides, then we may consider 
in place of $F$ the union of $F$ and all the $\e$-squares with centers $(x,y)$ in $Q\cap\e\Z^2$ and 
$$y\le \max \{z_2: z\in \Gamma\}.$$ 
The new set, pictured in Fig.~\ref{upper}, has both lower perimeter and less symmetric difference with $E^0_\e$. 

\begin{figure}[ht]
\centering
\def\svgwidth{280pt}
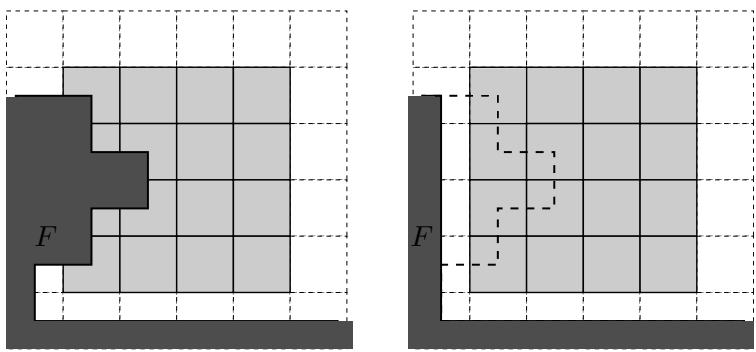 
\caption{Removal of $\partial F$ when intersecting one side.}\label{same}
\end{figure}

\begin{figure}[ht]
\centering
\def\svgwidth{280pt}
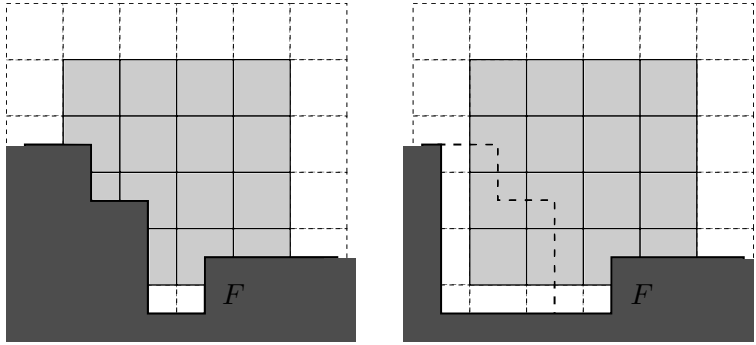  
\caption{Removal of $\partial F$ when intersecting two adjacent sides.}\label{adjacent}
\end{figure}

If a portion $\Gamma$ intersects $Q$ in exactly two points on the same side (horizontal or vertical) or adjacent sides,
then we may remove all the $\e$-squares with centers in the portion of $Q\cap F$ with boundary $\Gamma$.
The two cases are pictured in Fig.~\ref{same} and Fig.~\ref{adjacent}, respectively. 
This operation decreases the perimeter of at least $\e(\beta-\alpha)$, while at most increases the 
bulk term by ${1\over \tau}\e^3 N_\beta^2 c(L)$ ($c(L)$ given by (\ref{cielle})). The total change in the energy is thus
\begin{equation}\label{stima}
-\e(\beta-\alpha) +{1\over \gamma}\e^2N_\beta^2 c(L)\,,
\end{equation}
which is negative if $\e$ is small enough.
As a consequence, either $F\cap Q=\emptyset$ or $\partial F\cap Q$ is a horizontal segment.

The same type of analysis applies to the extremal squares, for which we deduce instead that $F\cap Q$ is a rectangle with one vertex coinciding with a vertex of $\widetilde R^\alpha$.

\begin{figure}[ht]
\centering
\def\svgwidth{420pt}
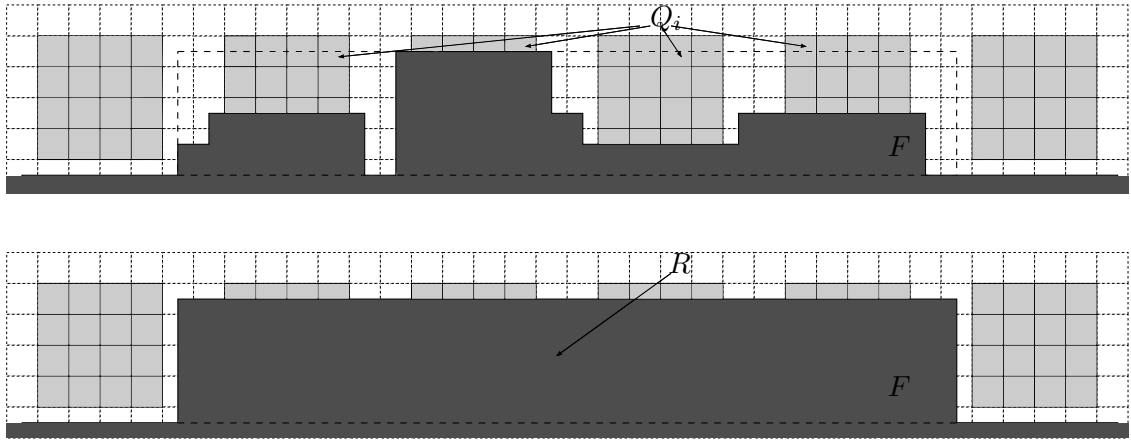  
\caption{Envelope of $\partial F$ in consecutive squares.}\label{more}
\end{figure}

\begin{figure}[ht]
\centering
\def\svgwidth{420pt}
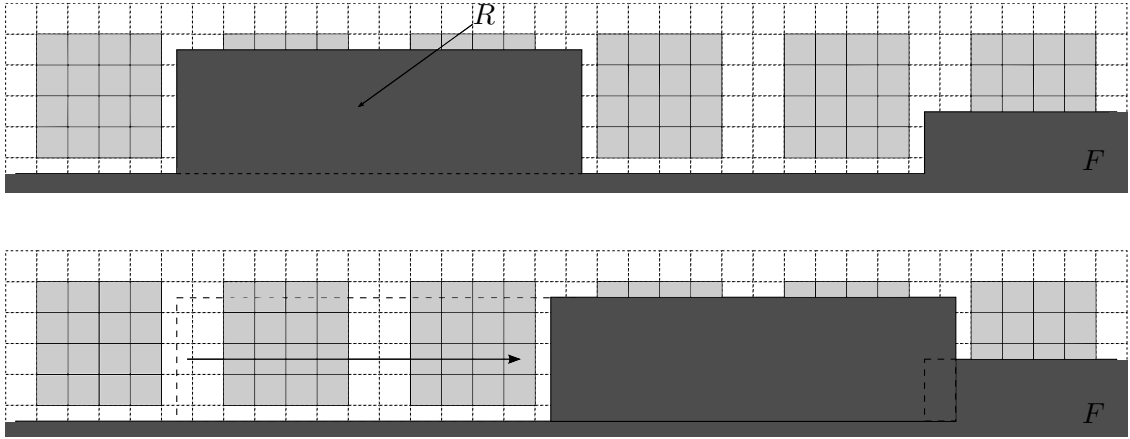 
\caption{Translation argument to join protrusions.}\label{right}
\end{figure}

\begin{figure}[ht]
\centering
\def\svgwidth{420pt}
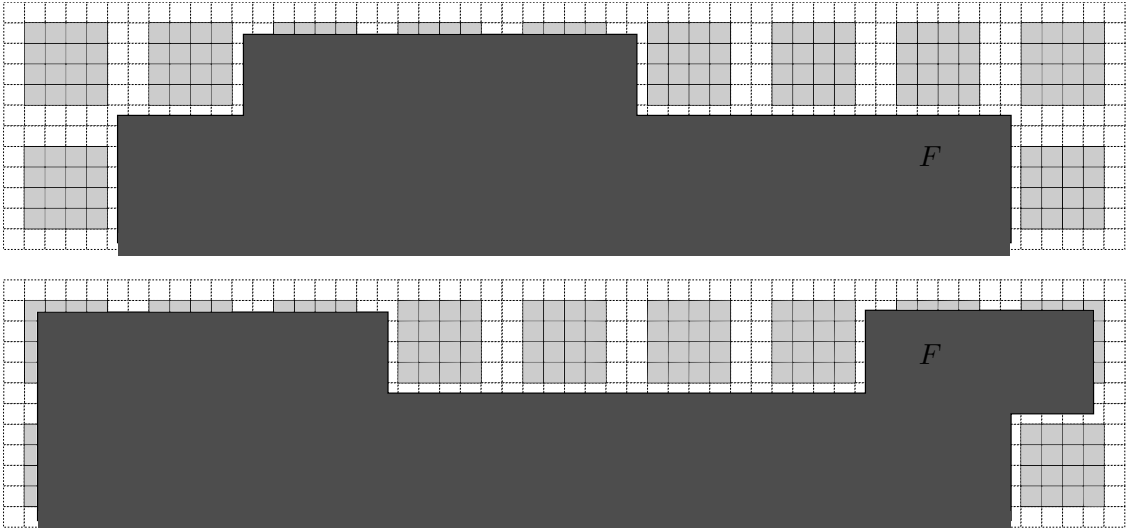 
\caption{Profiles of candidate minimal $F$.}\label{final}
\end{figure}

\bigskip We now consider the interaction of consecutive $\beta$-squares.
 Let $Q_1,\ldots, Q_K$ be a maximal array of consecutive $\beta$-squares with $F\cap Q_k\neq\emptyset$ for $k=1,\ldots, K$
and such that $Q_1$ is not a corner square.
If we substitute $F$ with $F\cup R$, where $R$ is the maximal rectangle of $\e$-squares containing all $F\cap Q_k$ and not intersecting other $\beta$-squares, then the corresponding energy has a not larger perimeter part, and a bulk part which is strictly lower if $F\cup R\neq F$. This substitution is pictured in Fig.~\ref{more}. 
If the subsequent $\beta$-squares $Q_{K+1}\ldots, Q_{K+K'}$  are a maximal array which do not intersect $F$, then we may further substitute $F\cup R$ with $(F\setminus R)\cup (R+\e \Nab K'(1,0))$, where we translate $R$ until it meets another portion of $F$ (if any).  This translation is pictured in Fig.~\ref{right}. 
Note that if it does meet another portion of $F$, then the change in energy is at most
\begin{equation}\label{stima2}
-2\e\alpha  +{1\over \gamma}\e^2N_\beta N_\alpha c(L)\,,
\end{equation}
which is negative if $\e$ is small enough. 
In this case at this point we may iterate this analysis since we now have a larger array of consecutive $\beta$-squares intersecting $F$. Note, moreover, that the same argument can be repeated shifting the rectangle $R$ to the left instead than to the right if energetically convenient. As a conclusion, we obtain that $F$ may only either intersect one array of consecutive $\beta$-squares, or two such arrays if they contain the two corner $\beta$-squares; i.e., we have one of the two situations pictured in Fig.~\ref{final}.

\smallskip
{\bf Step 4: all $\beta$-connections can be removed except those at the four corners.} At this point we are in the situations  pictured in Fig.~\ref{final}.
\noindent
If we are as in the upper figure, then by removing all $\e$-squares external to $R^\alpha$ the variation of the energy is less or equal than
\begin{equation*}
-(\beta-\alpha)(\Nb+1)N\epsilon+c(L)\frac{(N+1)N_{\alpha\beta}\Nb}{\gamma}\epsilon^2,
\end{equation*}
where $N$ is the number of modified $\beta$-squares. For $\epsilon$ small this variation is negative, showing that $F$ does not contain any protrusion.

\begin{figure}[ht]
\centering
\def\svgwidth{420pt}
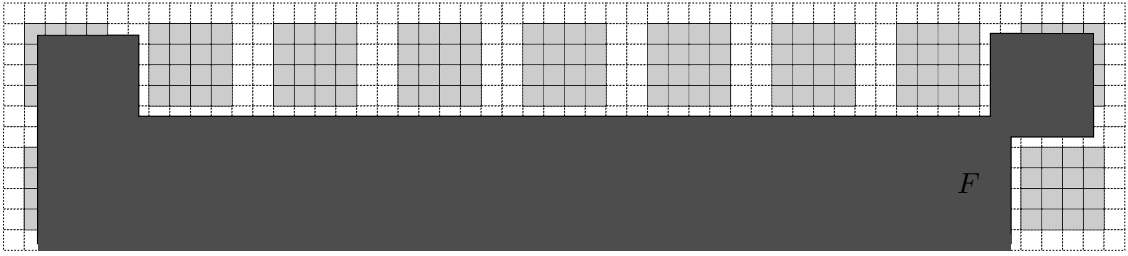
\caption{Removing $\beta$\hbox{-}connections except in the two $\beta$-squares at the corners.}\label{fig:12}
\end{figure}

\begin{figure}[ht]
\centering
\def\svgwidth{400pt}
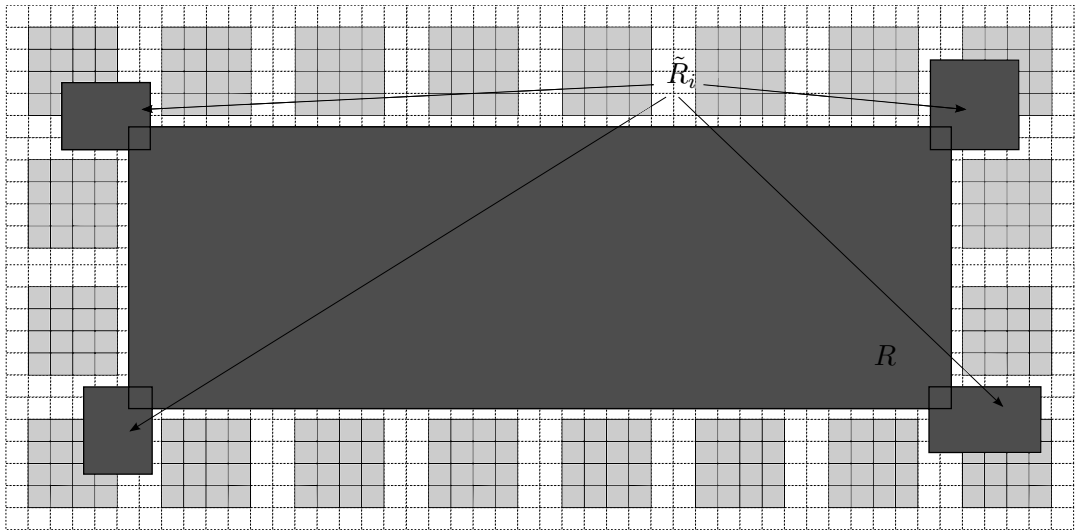
\caption{The set obtained in Step 4.}\label{fig:quasirectangle}
\end{figure}

If we are as in the lower figure, then we may remove all $\beta$\hbox{-}connections inside the border $\beta$-squares, except those in  the two periodicity squares at the corners as in Fig.~\ref{fig:12};
\noindent
the variation of the energy functional is less or equal than
\begin{equation*}
-(\beta-\alpha)(\Nb+1)N\epsilon+c(L)\frac{NN_{\alpha\beta}\Nb}{\gamma}\epsilon^2,
\end{equation*}
\\
where $N$ is the number of modified cells. For $\epsilon$ small this variation is negative, showing that the profile in Fig.~\ref{fig:12} is energetically convenient.
We can repeat this procedure for each side, and finally we obtain that $F$ is 
the union of a  coordinate $\alpha$-type rectangle $R$ and possibly one to four rectangles $\widetilde{R}_i, i=1,\dots,4$ of side lengths at most $\Nab\e$
such that  the intersection of $\widetilde{R}_i$ with each corner $\beta$-square is a rectangle (see Fig.~\ref{fig:quasirectangle}).

\smallskip
{\bf Step 5: conclusion.} It remains to prove that the rectangles $\widetilde{R}_i$ in the previous step are actually not there. 
This is immediately checked by comparing such an $F$ with $R^\alpha$: if $\widetilde{R}_i\neq \emptyset$ then by removing it the energy changes by at most by
$$
-2\beta\e+{1\over\gamma}c(L) \e^2\Nab^2,
$$
which is negative for small $\e$.

We finally note that all the estimates above can be iterated and hold uniformly as long as the sides of $E^k_\e$ are larger than $\delta$, since they depend only on $c(\delta)$.
\endproof

The proposition above shows that we may restrict our analysis to $\alpha$-type rectangles; indeed, for fixed $\e$ this assumption is not restrictive until the sides of the rectangles are larger than a constant, which vanishes as $\e\to 0$. 
As a consequence, once we suppose the convergence of the initial data, up to subsequences, the discrete 
motions $E_{\epsilon,\tau}(t)$ converge as $\e\to 0$ to a limit $E(t)$ such that $E$ is a rectangle for all $t$, up to its extinction time. Note, moreover, that it is not restrictive to suppose that also the initial data are  $\alpha$-type rectangles, up to substituting $E^0_\e$ with $E^1_\e$.
\begin{figure}[ht]
\centering
\def\svgwidth{300pt}
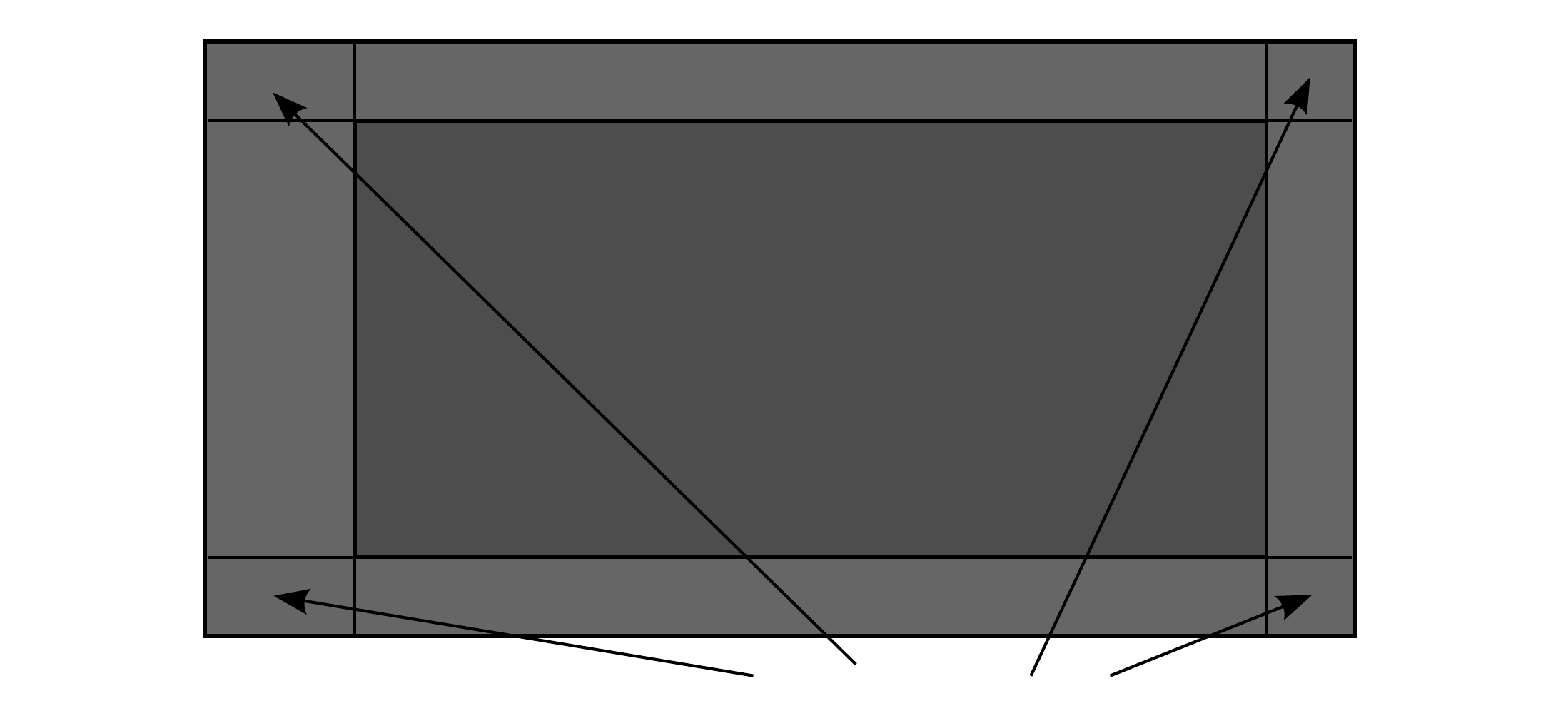
\caption{Picture of $E^{k+1}_\e$ inside $E^k_\e$.}\label{BNG-figure4}
\end{figure}

As shown in \cite{BGN}, the motion of each side of $E_\e^k$ can be studied separately, since the constraint of being an $\alpha$-type rectangle does not influence the argument therein, which consists in remarking that the bulk term due to the small corner rectangles in Fig.~\ref{BNG-figure4} is negligible.
As a consequence, we can describe the motion in terms of the length of the sides of $E_\e^k$. This will be done in the following sections.

\subsection{A new {pinning} threshold}\label{newpinning}
\begin{figure}[htbp]
\centering
\def\svgwidth{190pt}
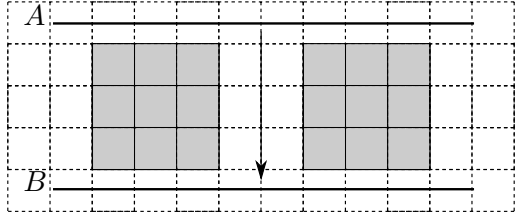
\caption{Motion is possible if the side can move at least by $(\Nb+1)\epsilon$.}
\label{fig:pinning}
\end{figure}

We first examine the case when the limit motion is trivial; i.e., all $E_k=E_\e^k$ are the same after a finite number of steps.
This will be done by computing the {\em pinning threshold}; i.e., the critical value of the side length $L$ above which it is energetically not favorable for a side to move. We recall that, in the case $\alpha=\beta$, this threshold is $$\widetilde{L}=2\alpha\gamma.$$ 
This value is obtained by computing the values for which a side of length $L$ may not move inwards of $\e$ by decreasing the energy. In our case, by the condition that $E_k$ be an $\alpha$-type rectangle, we have to impose instead that it is not energetically favorable to move inwards a side by $(\Nb+1)\varepsilon$ (see Fig.~\ref{fig:pinning}). 
We then  write the variation of the energy functional $\mathcal{F}^{\alpha,\beta}_{\epsilon,\tau}$ from configuration $A$ to configuration $B$ in Fig.~\ref{fig:pinning}, regarding a side of length $L$.
 If we impose it to be positive, we have

\begin{equation*}
-2(\Nb+1)\alpha\epsilon+\frac{1}{\tau}\displaystyle\sum_{k=1}^{\Nb+1}(k\epsilon)L\epsilon =
(\Nb+1)\epsilon\left[-2\alpha+\frac{L}{2\gamma}(\Nb+2)\right]\geq0
\end{equation*}
and we obtain the pinning threshold
\begin{equation}
\overline{L}:=\frac{4\gamma\alpha}{\Nb+2}.
\label{threshold}
\end{equation}
\\
Note that this threshold depends on $\Nb$ and not on the value $\beta>\alpha$ and that, if $\Nb=0$ (or, otherwise, $\alpha=\beta$), we recover the previous threshold $\widetilde L$.

\subsection{Definition of the effective velocity}\label{newvelocity}
As remarked above, up to an error vanishing as $\e\to0$, the motion of each side is independent of the other ones. 
As a consequence, its description can be reduced to a one-dimensional problem, where the unknown represents, e.g., the location of the left-hand vertical side of $E_k$.

Let $x_k$ represents the projection of this side of $E_k$ on the horizontal axis.
The location of $x_{k+1}$ depends on a minimization argument involving
$x_{k}$ and the length $L_k$ of the corresponding side of $E_{k}$.
However, we will see that this latter dependence is locally constant, except for a discrete set of values of $L_k$.
Indeed, for all $Y>0$ (which in our case will be of the form $Y=\gamma/L_k$), consider the minimum problems
\begin{equation}\label{minel}
\min\left\{-2\alpha N+\frac{N(N+1)}{2Y}:N\in\N, \quad [x+N]_{N_{\alpha\beta} }\in\Z_{\Na}\right\},
\end{equation}
for $x\in\{0,\ldots, \Nab\}$,
where $[z]_{N_{\alpha\beta}}$ denotes the congruence class of $z$ modulo $N_{\alpha\beta}$ and
\begin{equation*}
\Z_{\Na}=\left\{[0]_{N_{\alpha\beta}},\dots,[\Na-1]_{N_{\alpha\beta}}\right\}.
\end{equation*}
 Then the set of $Y>0$ for which (\ref{minel}) does not have a unique solution is discrete.
To check this it suffices to remark that the function to minimize 
$$-4\alpha XY+{X(X+1)}$$
is a parabola with vertex in 
$$
X={2\alpha Y}-{1\over 2}.
$$
The minimizers $N$ are points with $[x+N]_{N_{\alpha\beta} }\in\Z_{\Na}$ of minimal distance from the vertex $X$.
These are not unique in some cases: first if the vertex $X$ is equidistant from two consecutive points in ${\Z}_{\Na}$; i.e.,
if $${2\alpha  Y}-{1\over 2}\in {1\over 2}+\Z, $$ or, equivalently, 
\begin{equation}\label{evena}
Y\in{1\over 2\alpha}\Z.
\end{equation}

The second case is when we have two points in ${\Z}_{\Na}$ of minimal distance from $X$ which are not consecutive. In this case the distance between these points is $\Nb+1$, so that we have
$$
{2\alpha  Y}-{1\over 2}\in {\Nb+1\over 2}+\Z,
$$
or, equivalenly,
$$
Y\in{1\over 2\alpha}\Bigl({\Nb\over 2}+\Z\Bigr).
$$
If $\Nb$ is even then this condition is equivalent to (\ref{evena}), while if $\Nb$ is odd then we have 
\begin{equation}\label{odda}
Y\in {1\over 4\alpha}+{1\over 2\alpha}\Z.
\end{equation}

\def\Sb{S_\beta}

\begin{defn} We define the (possibly) {\em singular set} $S_{\beta}$ for problems (\ref{minel}) as 
$$
\Sb= {1\over 2\alpha}\Bigl(\Z\cup\Bigl({1\over 2}+\Z\Bigr)\Bigr).
$$
\end{defn}

We will examine the iterated minimizing scheme for $\gamma /L_k=\gamma /L\in (0,+\infty)\setminus S_\beta$ fixed, which reads
\begin{equation}
\begin{cases}
x_{k+1}^L=x_k^L+\overline{N}_k,&k\geq0\\
x_0^L=x^0
\end{cases}
\label{system}
\end{equation}
with $x^0\in\{0,1,\dots,N_{\alpha\beta}-1\}$ and $\overline{N}_k\in\N$ the minimizer of

\begin{equation}
\min\left\{-2\alpha N+\frac{1}{\gamma}\frac{N(N+1)}{2}L:N\in\N, \quad [x^L_k+N]_{N_{\alpha\beta} }\in\Z_{\Na}\right\},
\label{scheme}
\end{equation}
which is unique by the requirement that $\gamma/L\not\in S_\beta$.

After at most $\Na$ steps, $\{x_{k}^L\}_{k\geq0}$ is \emph{periodic modulo} $N_{\alpha\beta}$, as expressed by the following proposition.

\begin{prop}\label{perino}
There exist integers $\overline{k}\le \Na, M\leq \Na$ and $n\geq1$ such that
\begin{equation}\label{periodic}
x_{k+M}^L=x_{k}^L+n\,N_{\alpha\beta}
\qquad\hbox{ for all } k\geq\overline{k}.
\end{equation}
Moreover, the quotient $M/n$ depends only on $\gamma/L$.

\end{prop}

\proof

First remark that, if $x_k^L$ is defined recursively by (\ref{system}), we have
\begin{equation*}
 [x_k^L]_{N_{\alpha\beta}}\in{\Z_{\Na}}\qquad\hbox{ for all }k\ge 1.
\end{equation*}

Since $\#\Z_{\Na}=\Na$, there exist integers $0\leq j\leq \Na$ and $l>j$, with $l-j\leq \Na$, such that

\begin{equation}
[x^L_j]_{N_{\alpha\beta}}=[x^L_{l}]_{N_{\alpha\beta}}.
\label{equality}
\end{equation}
Let $l$ be the minimal such $l$.
Define $\overline{k}=j$, $M=l-j$ and $n=\displaystyle\frac{x^L_l-x^L_j}{N_{\alpha\beta}}$ to obtain (\ref{periodic}).

\smallskip
It remains to show the last statement of the theorem. It suffices to show that the quotient is independent of $x_0$.
We start by proving a monotonicity property of the orbits defined in (\ref{system}) 
with respect to the initial datum: if $\{x_k\}$ and  $\{x'_k\}$ are orbits obtained as above, we have
\begin{equation}
\text{if } x_0\le x'_0, \text{then }x_k\leq x'_k,\quad \hbox{for all }k\geq1.
\label{monotone}
\end{equation}
This can be seen iteratively from (\ref{minel}) since the problems with $x=x_{k-1}$ and $x=x'_{k-1}$ 
consist in a constrained minimization of a parabola and its translation by $x'_{k-1}-x_{k-1}$, and, as
previously remarks, the minimizer 
in (\ref{minel}) is the closest point to the  vertex of the parabola with $[x+N]_{N_{\alpha\beta} }\in\Z_{\Na}$.

Consider the orbits with initial data $x_0$, $x'_0$ and $x_0+\Nab$,
and let $n(x)$ and $M(x)$ denote the indices above with initial datum $x\in\{x_0,x'_0,x_0+\Nab\}$.
Since the orbit with initial datum $x_0+\Nab$ is the translation by $\Nab$ of the one with initial datum $x_0$,
we have $n(x_0+\Nab)=n(x_0)$ and $M(x_0+\Nab)=M(x_0)$. Taking into account the ordering of the initial conditions
$$
x_0\le x'_0\le x_0+\Nab,
$$
by (\ref{periodic}) for $k_0$ sufficiently large and taking $k=k_0+T M(x_0) M(x'_0)$ with $T\in\N$,
from $x_k\leq x'_k\leq x_k+\Nab$
we get
\begin{eqnarray*}
x_{k_0}+Tn(x_0)M(x'_0)\Nab\Nab&\le&  x'_{k_0}+Tn(x'_0)M(x_0)\Nab\\
&\le& x_{k_0}+Tn(x_0)M(x'_0)\Nab+\Nab.
\end{eqnarray*}
In order that this inequality hold for all $T\ge 1$ we must have 
$$
n(x_0)M(x'_0)=n(x'_0)M(x_0),
$$
which is the desired equality.
\endproof

\begin{defn}[Effective velocity]\label{velofo} \rm
We define the 
\emph{effective velocity function}
$f:(0,+\infty)\setminus S_\beta\longrightarrow[0,+\infty)$ by setting 
\begin{equation}
f(Y)=\frac{nN_{\alpha\beta}}{M},
\label{velocityfunction}
\end{equation}
with $M$ and $n$ in (\ref{periodic}) defined by $L$ and $\gamma$ such that $Y=\gamma/L$. By Proposition \ref{perino} this is a good definition.
\label{effvel}
\end{defn}

\begin{oss}\label{meva}
The terminology for formula
 (\ref{velocityfunction}) is motivated by the fact that we can define the velocity of a side as a \emph{mean velocity}
 averaging on a period; that is,
\begin{equation}
v=\frac{nN_{\alpha\beta}\epsilon}{M\tau}.
\label{meanvel}
\end{equation}
In (\ref{meanvel}) the velocity is the ratio between the minimal (periodic) displacement of the side and the product of the time\hbox{-}scale $\tau$ and the number of steps necessary to describe the minimal period, each of which considered as a 1\hbox{-}time step. 
\end{oss}

\begin{oss}[Properties of the velocity function $f$] The velocity function $f$ has the following properties:
\begin{description}
\item[(a)] $f$ is constant on each interval contained in its domain;
\item[(b)] $f(Y)=0$ if $$
Y<\overline Y:=\frac{\Nb+2}{4\alpha};
$$
in particular
$$
\lim_{\gamma\to 0^+} {1\over \gamma} f\Bigl({\gamma\over L}\Bigr)=0\,.
$$
Note that $(0,\overline Y)\cap S_\beta\neq\emptyset$;
\item[{(}c)] $f(Y)$ is a rational value;
\item[(d)] $f$ is non decreasing;
\item[(e)] we have
$$
\lim_{\gamma\to+\infty} {1\over \gamma} f\Bigl({\gamma\over L}\Bigr)={2\alpha\over L}\,.
$$
\item[(f)] $f(Y)$ is independent of $\beta$ but depends on $\Nb$.
\end{description}

{\bf (a)} holds since on each component of $(0,+\infty)\setminus S_\beta$ the minimum problems (\ref{minel}) have a unique solution independent of $Y$, so that the values $n$ and $M$ in Proposition \ref{perino} are independent of $Y$. Note, however, that $f(Y)$ may be equal on neighboring components since the corresponding $n$ and $M$ may be equal even without uniqueness in (\ref{perino});

{\bf (b)} holds since we have $\overline{Y}=\gamma/\overline L$,  where $\overline{L}$ is the pinning threshold (\ref{threshold}),
and the computation of the pinning threshold is equivalent to the requirement that the orbit be constant after a finite number of steps;

{\bf (c)} is immediate from the formula for $f(Y)$;

{\bf (d)} is again a consequence of the fact that (\ref{minel}) are minimum problems related to a parabola with vertex in $2{\alpha Y}-{1\over 2}$ and the latter is an increasing function of $Y$;

{\bf (e)} using the same argument as in {\bf(d)} above, we deduce in particular that
$$
\Bigl|\overline {N}_k-{2\alpha Y}+{1\over 2}\Bigr|\le \Nb,
$$
which for $Y=\gamma/L$ implies that 
$$
{2\alpha\over L}-{2\Nb+1\over 2\gamma}\le {1\over \gamma} f\Bigl({\gamma\over L}\Bigr)\le {2\alpha\over L}+{2\Nb+1\over 2\gamma},
$$
and the desired equality letting $\gamma\to+\infty$;

{\bf (f)} is an immediate consequence of the definition of $f(Y)$.
\end{oss}

\begin{oss}\label{compar}
Let $\gamma/L\in S_\beta$, and let $\{x^L_k\}$ be defined by (\ref{system}) with $\overline {N}_k$ chosen to be a minimizer of (\ref{scheme}), which may be not unique. Then arguing by monotonicity as in {\bf(d)} above, we have $ x^{L^+}_k\le x^L_k\le  x^{L^-}_k$, where $L^\pm$ are any two values with $L^-<L<L^+$ and $\gamma/L^\pm$ belonging to the two intervals of $(0,+\infty)\setminus S_\beta$ with one endpoint equal to $L$, and $\{x^{L^\pm}_k\}$ have the same initial data.
\end{oss}

\subsection{Description of the homogenized limit motion.}\label{limitmotion}

The following characterization of any limit motion holds. 

\begin{thm}\label{limitmotion1}For all $\epsilon>0$, let $E^0_\epsilon\in\D_\epsilon$ be a coordinate rectangle with sides $S^0_{1,\epsilon},\dots,S^0_{4,\epsilon}$. Assume also that
\begin{equation*}
\lim_{\epsilon\to0^+}\emph{d}_\H(E^0_\epsilon,E_0)=0
\end{equation*}
for some fixed coordinate rectangle $E_0$. Let $\gamma>0$ be fixed and let  $E_\e(t)= E_{\e,\gamma\e}(t)$ be the piecewise-constant motion with initial datum $E^0_\e$ defined in {\rm(\ref{disefo})}.
Then, up to a subsequence, $E_\epsilon(t)$ converges as $\epsilon\to0$ to $E(t)$, where $E(t)$ is a coordinate rectangle with sides $S_i(t)$ and such that $E(0)=E_0$. Each $S_i$ moves inward with velocity $v_i(t)$ satisfying
\begin{equation}
v_i(t)
\in\left[\displaystyle\frac{1}{\gamma}f\biggl({\gamma\over L_i(t)}\biggr)^-,\displaystyle\frac{1}{\gamma}f\biggl({\gamma\over L_i(t)}\biggr)^+\right],
\label{vl}
\end{equation}
where $f$ is given by Definition {\rm\ref{velofo}}, $L_i(t):=\H^1(S_i(t))$ denotes the length of the side $S_i(t)$, until the extinction time when $L_i(t)=0$, and $f(Y)^-,f(Y)^+$ are the upper and lower limits of the effective-velocity function at $Y\in (0,+\infty)$.
\end{thm}
\proof  We will apply the results of the previous sections with $\tau=\gamma\e$.
Let $S_{\e,i}(t)$ be the sides of $E_\e(t)$, and let $L_{i,\e}^k= \H^1(S_{\e,i}(k\tau))$; i.e., $L_{i,\e}^k$ is the length of the $i$-th side of $E^k_\e$ in the notation of the previous sections. If $\Delta S^k_{\e,i}= \text{d}_\H(S_{\e,i}(\gamma\e k),S_{\e,i}(\gamma\e (k+1))$ denotes the distance from corresponding sides of $E^k_\e$ then note that
$$
L_{i,\e}^{k+1}-L_{i,\e}^k= -\bigl(\Delta S^k_{\e,i-1}+ \Delta S^k_{\e,i+1}\bigr)
$$
(where the indices $i$ rotate cyclically).
By (\ref{cielle}) we have 
$$
{\Delta S^k_{\e,i}\over \tau}\le {c_1\over L_{i,\e}^{k}}+ c_2.
$$
This implies that if we define $L_{i,\epsilon}(t)$ as the affine interpolation 
in $[k\tau,(k+1)\tau]$ of the values $L_{i,\epsilon}^k$, then $L_{i,\epsilon}(t)$ is a decreasing continuous function of $t$ and the sequence is uniformly Lipschitz continuous on all intervals $[0,T]$ such that $L_{i,\epsilon}(T)\geq c>0$. Hence it converges (up to a subsequence), as $\epsilon\to0$, to a function $L_i(t)$, which is also decreasing. It follows that $E_\epsilon(t)$ converges as $\epsilon\to0$, up to a subsequence and in the Hausdorff sense, to a limit rectangle $E(t)$, for all $t\geq0$.

It remains to justify formula (\ref{vl}) for the velocity $v_i$ of the side $S_i(t)$. Let $[t^-, t^+]$ and $L^\pm _i$ be such that
$\gamma/L^\pm_i\in (0,+\infty)\setminus S_\beta$ and 
$$
L^-_i < L_i(t)< L^+_i \qquad \hbox{ for } t^-\!\le t\le t^+.
$$
Then the corresponding $L_{i,\e}(t)$ satisfy the same inequalities for $\e$ small enough. By Remarks \ref{compar} and \ref{meva} we then have
$$
{1\over \gamma} f\Bigl({\gamma\over L^+}\Bigr)\le v_i(t)\le {1\over \gamma} f\Bigl({\gamma\over L^-}\Bigr)\quad\hbox{ for } t^-\!\le t\le t^+.
$$
By optimizing in $L^\pm$, and recalling that $f$ is not decreasing, we obtain (\ref{vl}).
\endproof

\begin{thm}[Unique limit motions]
Let $E_\epsilon,E_0$ be as in the statement of Theorem~{\rm \ref{limitmotion}}. Assume in addition that the lengths $L^0_1,L^0_2$ of the sides of the initial set $E_0$ satisfy one of the following conditions (we assume that $L^0_1\leq L^0_2$):
\begin{itemize}
\item[\emph{(a)}] $L^0_1,L^0_2>\displaystyle\frac{4\alpha\gamma}{\Nb+2}$ \emph{(total pinning)};
\item[\emph{(b)}] $L^0_1<\displaystyle\frac{4\alpha\gamma}{\Nb+2}$ and $L^0_2\leq\displaystyle\frac{4\alpha\gamma}{\Nb+2}$ \emph{(vanishing in finite time)};
\end{itemize}
then $E_\epsilon(t)$ converges locally in time to $E(t)$ as $\epsilon\to 0$, where $E(t)$ is the unique rectangle with sides of lengths $L_1(t)$ and $L_2(t)$ which solve the following system of ordinary differential equations

\begin{equation}\label{unita}
\begin{cases}\displaystyle
\dot{L}_1(t)=-{2\over \gamma}\,f\Bigl({\gamma\over L_2(t)}\Bigr)\\
\\ \displaystyle
\dot{L}_2(t)=-{2\over \gamma}\, f\Bigl({\gamma\over L_1(t)}\Bigr)
\end{cases}
\end{equation}
for almost every $t$, with initial conditions $L_1(0)=L^0_1$ and $L_2(0)=L^0_2$, where $f$ is given by Definition {\rm\ref{velofo}}.
\end{thm}

\proof
In case (a) the statement follows by Theorem~\ref{limitmotion1} noticing that we have $v_1(t)=v_2(t)=0$ for all $t\geq0$, which is equivalent to $\dot{L}_1=\dot{L}_2=0$. 

In case (b) the lengths of $L_i$ are strictly decreasing until the extinction time. This implies that the set of $t$ such that
$f ({\gamma/ L_i(t)} )^-\neq f ({\gamma/ L_i(t)} )^+$ is negligible, and (\ref{unita}) follows since  $\dot L_i=-2 v_{i+1}$.
\endproof

\begin{oss}[general evolutions]
More general initial data can be considered. Since their treatment follows from Theorem \ref{limitmotion1} as in \cite{BGN}, we do not include the details. We only recall that:

$\bullet$ all velocities $v_i$ satisfying (\ref{vl}) can be obtained, with a proper choice of the initial data $E_{0,\e}$;

$\bullet$ if we take initial data $E_0$ coordinate polyrectangles then the motion can be characterized with the same velocities, with the convention that convex sides move inwards, concave sides move outwards, other sides remain pinned;

$\bullet$ more general  initial data $E_0$ can be dealt with once we remark that at level $\e$ the assumption that $E_{0,\e}$ is a polyrectangle is always satisfied.
\end{oss}

\section{Computation of the velocity function}
The velocity function in Definition \ref{effvel} may be not easily described 
 for generic $\Na$ and $\Nb$. In this section we compute it, by means of algebraic formulas, in the simpler cases $\Nb=1$ and $\Nb=2$, with varying $\Na$. These are prototypes for the cases $\Nb$ odd and $\Nb$ even, respectively. We also give two easy examples for $\Na$ fixed and equal to $1$, and we compare the new velocity function with the homogeneous case showing that the inhomogeneities in the lattice may accelerate or decelerate the motion. We can assume, without loss of generality, that $\gamma=1$.
 
\subsection{The case $\Nb=1$.}\label{N1}
Let $Y>\overline{Y}=\frac{3}{4\alpha}$. We assume also that $Y$ is not in the singular set; 
i.e.,
\begin{equation*}
Y\not\in\left\{\frac{k+j(\Na+1)}{2\alpha},k=1,\dots,\Na-1,j\ge 0\right\}\cup\left\{\frac{\Na+(2j+1)(\Na+1)}{4\alpha},j\ge 0\right\}.
\end{equation*}

As shown by Proposition~\ref{perino}, the minimal period is independent of the starting point of the orbits, so there is no restriction to assume that $x^0=0$ in (\ref{system})-(\ref{scheme}). 
We divide the analysis in the three cases (a), (b) and {(}c) below.

\bigskip
$(a)$ If $Y\in \left(\displaystyle\frac{k+j(\Na+1)}{2\alpha},\frac{k+1+j(\Na+1)}{2\alpha}\right), k=1,2,\dots,\Na-1,j\ge 0$, then we denote the minimizer of problem (\ref{scheme}) in the homogeneous case $\Nb=0$ by $N=k+j(\Na+1)$. The velocity function $f(Y)$ will be characterized by algebraic relations between $N$ and $\Na$. 
We have two sub-cases:

$(a_1)$ $N$ and $\Na+1$ are coprime. In this case, by iterating the scheme (\ref{scheme}), after at most $\Na$ steps the side encounters a defect, that is

\begin{equation*}
[nN]_{\Na+1}=[\Na]_{\Na+1}
\end{equation*}
for some $1\leq n\leq\Na$.
In this case, we denote by $\bar{n}\ge0$ the minimal solution of the congruence equation
\begin{equation}
nN\equiv\Na\quad \text{mod $(\Na+1),$   $n\geq1$},
\label{congr0}
\end{equation}
and $\bar{k}\geq0$ is given by
\begin{equation*}
\bar{k}=\frac{\bar{n}N-\Na}{\Na+1}.
\end{equation*}

If $Y\in\left(\displaystyle\frac{k+j(\Na+1)}{2\alpha},\frac{2k+2j(\Na+1)+1}{4\alpha}\right)$, then the location of the side at step $n$ is at $\Na-1+\bar{k}(\Na+1)$ (which is equal to $-2$ modulo $\Na+1$). 

This computation shows that we can limit our analysis to periodic orbits modulo $\Na+1$ with initial datum equal to $-2$ (or, equivalently, $\Na-1$). The period of such orbits is obtained as follows. We solve the congruence equation

\begin{equation}
nN\equiv1\quad \text{mod }\Na+1, 
\label{congr1}
\end{equation}
for $n\geq1$ and denote by $n_{\text{min}}$ the minimal positive solution of equation (\ref{congr1}); that is, the minimal positive integer in the class
\begin{equation*}
\left[N^{\varphi(\Na+1)-1}\right]_{\text{mod $\Na+1$}}.
\end{equation*}
The function $\varphi(n)$ is \emph{Euler's totient function} and it counts the integers $m$ such that $1\leq m<n$ and $m$ has no common divisors with $n$.
If we define
\begin{equation*}
k_{\text{min}}=\frac{n_{\text{min}}N-1}{\Na+1},
\end{equation*}
then we have that
\begin{equation}
\begin{split}
f(Y)&=\frac{k_{\text{min}}(\Na+1)}{\displaystyle\frac{k_{\text{min}}(\Na+1)+1}{\lfloor2\alpha Y\rfloor}}=\left(\frac{k_{\text{min}}(\Na+1)}{k_{\text{min}}(\Na+1)+1}\right){\lfloor2\alpha Y\rfloor}\\
&=\left(\frac{1}{1+\frac{1}{k_{\text{min}}(\Na+1)}}\right){\lfloor2\alpha Y\rfloor}.
\end{split}
\label{dec}
\end{equation}
Note that $f(Y)<\lfloor2\alpha Y\rfloor$, so that the velocity of the side reduces (deceleration) with respect to the homogeneous case. \\

Suppose now that $Y\in\left(\displaystyle\frac{2k+2j(\Na+1)+1}{4\alpha},\frac{k+1+j(\Na+1)}{2\alpha}\right)$, then the location of the side at step $n$ is $\Na+1+\bar{k}(\Na+1)$, which is equal to $0$ modulo $\Na+1$.
We have that
\begin{equation}
f(Y)=\left(\frac{(\bar{k}+1)(\Na+1)}{(\bar{k}+1)(\Na+1)-1}\right){\lfloor2\alpha Y\rfloor}=\left(\frac{1}{1-\frac{1}{(\bar{k}+1)(\Na+1)}}\right){\lfloor2\alpha Y\rfloor}.
\end{equation}
Note that $f(Y)>{\lfloor2\alpha Y\rfloor}$, so the velocity of the side increases (acceleration) with respect to the homogeneous case.

$(a_2)$ $N$ and $\Na+1$ are not coprime. In this case the side does not meet any $\beta$\hbox{-}bond and the velocity function has the same value as in the homogeneous case, i.e.
\begin{equation*}
f(Y)=\lfloor2\alpha Y\rfloor.
\end{equation*}

\medskip
$(b)$ If $Y\in\left(\displaystyle\frac{\Na+j(\Na+1)}{2\alpha},\frac{\Na+(2j+1)(\Na+1)}{4\alpha}\right)$ then  we argue as in $(a_1)$.\\

\smallskip
$(c)$ If $Y\in\left(\displaystyle\frac{\Na+(2j+1)(\Na+1)}{4\alpha},\frac{1+(j+1)(\Na+1)}{2\alpha}\right)$, then

\begin{equation*}
f(Y)=\Na+1+j(\Na+1).
\end{equation*}
\\
Note that $f(Y)>\lfloor 2\alpha Y\rfloor$ if $Y\in\left(\displaystyle\frac{\Na+(2j+1)(\Na+1)}{4\alpha},\frac{(j+1)(\Na+1)}{2\alpha}\right)$, while $f(Y)=\lfloor 2\alpha Y\rfloor$ if $Y\in\left(\displaystyle\frac{(j+1)(\Na+1)}{2\alpha},\frac{1+(j+1)(\Na+1)}{2\alpha}\right)$.
\begin{example}[The case $\Na=\Nb=1$]
In this case the velocity function is given by

\begin{equation*}
f(Y)=
\begin{cases}
0&\mbox{if}\quad Y<\displaystyle\frac{3}{4\alpha},\\
\\
\displaystyle{2k}&\mbox{if}\quad Y\in\left(\displaystyle\frac{4k-1}{4\alpha},\displaystyle\frac{4k+3}{4\alpha}\right),\quad k\geq1;
\end{cases}
\end{equation*}
i.e., 
$$
f(Y)= 2\Bigl\lfloor \alpha Y+{1\over 4}\Bigr\rfloor.
$$
\end{example}

\subsection{The case $\Nb=2$}\label{N2}

We now study the case $\Nb=2$. Let $Y>\overline{Y}=\displaystyle\frac{1}{\alpha}$ and we assume also that $Y$ is not in the singular set, i.e.,

\begin{equation*}
Y\not\in\left\{\frac{k+j(\Na+2)}{2\alpha},k=1,\dots,\Na-1,j\ge 0\right\}\cup\left\{\frac{\Na+1+j(\Na+2)}{2\alpha},j\ge 0\right\}.
\end{equation*}

\smallskip
$(a)$ If $Y\in \left(\displaystyle\frac{k+j(\Na+2)}{2\alpha},\frac{k+1+j(\Na+2)}{2\alpha}\right), k=1,2,\dots,\Na-2,j\ge 0$, then $N=k+j(\Na+2)$ and we have two sub-cases:

$(a_1)$ $N$ and $\Na+2$ are coprime. We compute $\bar{k}=\min(k_1,k_2)\geq0$, where $k_1$ is the minimal positive solution of the congruence equation 
\begin{equation*}
kN\equiv\Na\quad \text{mod }\Na+2, 
\end{equation*}
and $k_2$ is the minimal positive solution of the congruence equation 

\begin{equation*}
kN\equiv \Na+1\quad \text{mod }\Na+2;
\end{equation*}
that is $k_1$ is the minimal positive integer in the class $\left[\Na N^{\varphi(\Na+2)-1}\right]_{\text{mod $\Na+2$}}$ and $k_2$ is the minimal positive integer in the class $\left[(\Na+1)N^{\varphi(\Na+2)-1}\right]_{\text{mod $\Na+2$}}$.

If $\bar{k}=k_1$, then 
\begin{equation}
f(Y)=\left(\frac{k_1(\Na+2)}{k_1(\Na+2)+1}\right){\lfloor 2\alpha Y\rfloor}=\left(\frac{1}{1+\frac{1}{k_1(\Na+2)}}\right){\lfloor 2\alpha Y\rfloor},
\end{equation}
and $f(Y)<\lfloor 2\alpha Y\rfloor$.

If $\bar{k}=k_2$, then
\begin{equation}
f(Y)=\left(\frac{k_2(\Na+2)}{k_2(\Na+2)-1}\right){\lfloor 2\alpha Y\rfloor}=\left(\frac{1}{1-\frac{1}{k_2(\Na+2)}}\right){\lfloor 2\alpha Y\rfloor},
\end{equation}
and $f(Y)>\lfloor 2\alpha Y\rfloor$.

$(a_2)$ $N$ and $\Na+2$ are not coprime. In this case
\begin{equation*}
f(Y)=\lfloor 2\alpha Y\rfloor
\end{equation*}
as in the homogeneous case.

\smallskip
$(b)$ If $Y\in \left(\displaystyle\frac{\Na+1+j(\Na+2)}{2\alpha},\frac{(j+1)(\Na+2)+1}{2\alpha}\right),j\ge 0$, then
\begin{equation*}
f(Y)=(j+1)(\Na+2).
\end{equation*}

Note that, in this case, if $Y\in \left(\displaystyle\frac{\Na+1+j(\Na+2)}{2\alpha},\frac{(j+1)(\Na+2)}{2\alpha}\right)$ then $f(Y)>\lfloor 2\alpha Y\rfloor$, while if 
$Y\in \left(\displaystyle\frac{(j+1)(\Na+2)}{2\alpha},\frac{(j+1)(\Na+2)+1}{2\alpha}\right)$ then $f(Y)=\lfloor 2\alpha Y\rfloor$.

\smallskip
$(c)$ If $Y\in \left(\displaystyle\frac{\Na-1+j(\Na+2)}{2\alpha},\frac{\Na+1+j(\Na+2)}{2\alpha}\right),j\ge 0$, then we may argue as in case $(a)$.

\begin{example}[The case $\Na=1, \Nb=2$]
The velocity function is given by
\begin{equation*}
f(Y)=
\begin{cases}
0&\mbox{if}\quad Y<\displaystyle\frac{1}{\alpha},\\
\\
\displaystyle{3k}&\mbox{if}\quad Y\in\left(\displaystyle\frac{3k-1}{2\alpha},\displaystyle\frac{3k+2}{2\alpha}\right),\quad k\geq1;
\end{cases}
\end{equation*}
i.e., 
$$
f(Y)=
3\Bigl\lfloor {2\over 3}\alpha Y+{1\over 3}\Bigr\rfloor.
$$
\end{example}

\end{document}